\documentclass[12pt]{article}
\usepackage{geometry}
\usepackage{graphicx}
\usepackage[dvips,arrow,matrix,ps,color,line]{xy}
\usepackage{theorem,amsfonts,amssymb,amscd,graphics,shapepar,bm}
\usepackage{epstopdf}
\usepackage{hyperref}

\hypersetup{ pdftitle={Catalan},
      pdfauthor={Taise Santiago},
      pdfsubject={Integral on Grassmanninans},
    colorlinks = true,
    linkcolor = blue,
    anchorcolor = red,
    citecolor = blue,
    filecolor = red,
    pagecolor = red,
    urlcolor = blue
}

\makeindex

\def\NN{\mathbb N}
\def\ZZ{\mathbb Z}

\def\PP{\mathbb P}
\def\ep{{\epsilon}}
\def\w{\wedge }
\def\lra{\longrightarrow}
\def\proof{\noindent{\bf Proof.}\,\,}
\def\qed{{\hfill\vrule height4pt width4pt depth0pt}\medskip}

\theoremstyle{change} \theorembodyfont{\rmfamily}

\newtheorem{claim}{}[section]

\title{``Catalan Traffic" and Integrals on the Grassmannian of Lines\thanks{AMS 2000 Math. Subject
Classification: 14M15, 14N15, 05A15, 05A19.}}
\author{Ta\'{i}se Santiago Costa Oliveira \thanks{This work was supported in part by ScuDo - Politecnico di Torino, FAPESB proc. n° 8057/2006 and CNPq proc. n° 350259/2006-2.}}
\date{}                                           

\begin{document}
\maketitle

\begin{abstract}
\noindent We prove that certain  numbers occurring in a problem of
paths enumeration, studied by Niederhausen in~\cite{Nied} (see
also~\cite{Stan1}), are top intersection numbers in the cohomology
ring of the grassmannian of the lines in the complex projective
$(n+1)$-space.

\end{abstract}

\section{Introduction}
\claim{} The {\em Catalan's numbers}
\[
C_n={1\over n+1}{2n\choose n}, \ \ \textnormal{for all}  \ \
n\in\NN
\]
occur in several combinatorial situations (see e.g.~\cite{Stan}),
in particular in {\em lattice path enumeration}. It is well known,
for instance, that  $C_n$ is the number of lattice paths contained
in $S:=\{(m,n)\in\ZZ\times \ZZ\,|\, 0\leq m\leq n\}$ from $(0,0)$
to $(n,n)\in S$, allowing unitary steps  only, along the ``horizontal"  or
``vertical"  directions.

Within this context, the aim of this paper is to make some remarks
on the occurrence of Catalan's numbers in a {\em traffic game}
(``Catalan traffic at the beach") constructed by
Niederhausen~\cite{Nied}. One is given of a {\em city map} ${\cal
C}$  (a lattice in $\ZZ^2$) with some {\em gates} and {\em road
blocks} (null traffic points). The traffic rules are as follows.
First, no path can cross and go beyond the {\em beach} (the line
$m-n=0$ in the   $(m,n)$ $\ZZ$-plane). Furthermore:
\begin{enumerate}
    \item At lattice points strictly ``below" the line $2m+n=0$, only
    North ($\uparrow$) or West ($\leftarrow$) directions are allowed;
    \item At lattice points strictly ``above" the line $2m+n=0$,
  only  East ($\rightarrow$) or NE  ($\nearrow$) directions are allowed;
    \item All the points $(m,n)\in\ZZ^2$ lying on the line $2m+n=1$ are inaccessible (road blocks $\blacksquare$).
    \item On the line $2m+n=0$ (gates), allow $W$($\leftarrow$), $E$ ($\rightarrow$), and
    $NE$ ( $\nearrow$) (because of the road blocks at $2m+n=1$) (see~\cite{Nied}, p.~2)
\end{enumerate}

The diagram of the city map is depicted below: it is the same as
in~\cite{Nied} after a harmless counterclockwise rotation of 90
degrees.

\begin{center}
 \includegraphics[width=9cm]{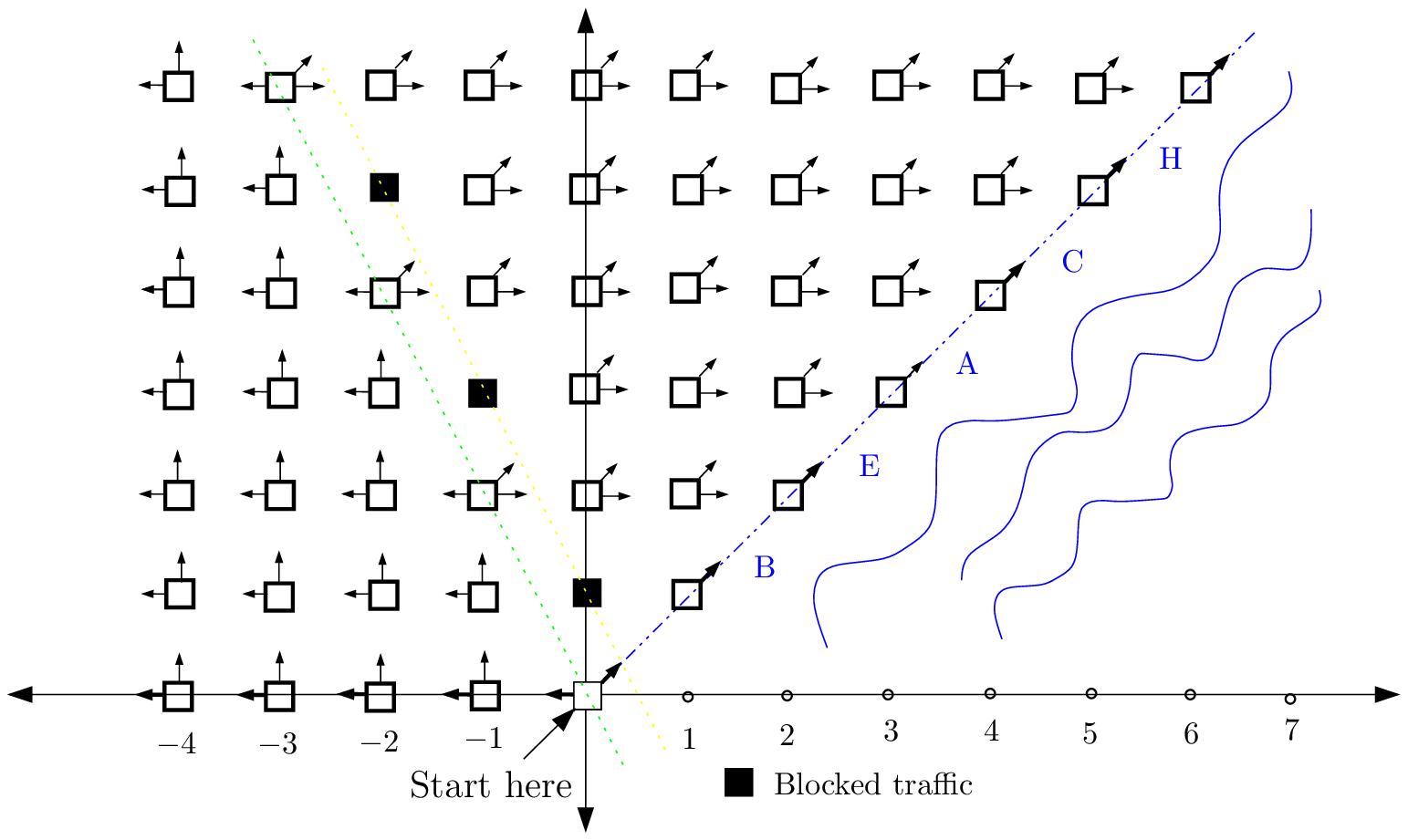}\\
 \medskip
 Fig.~1. City Map.
\end{center}

The problem, solved by Niederhausen, consists in finding the
number of all distinct paths joining the origin to any point in
the domain of ${\cal C}$, compatibly with the constraints. Attaching to each point of the lattice
the number of such paths one gets the following diagram:

\begin{center}
 \includegraphics[width=9cm]{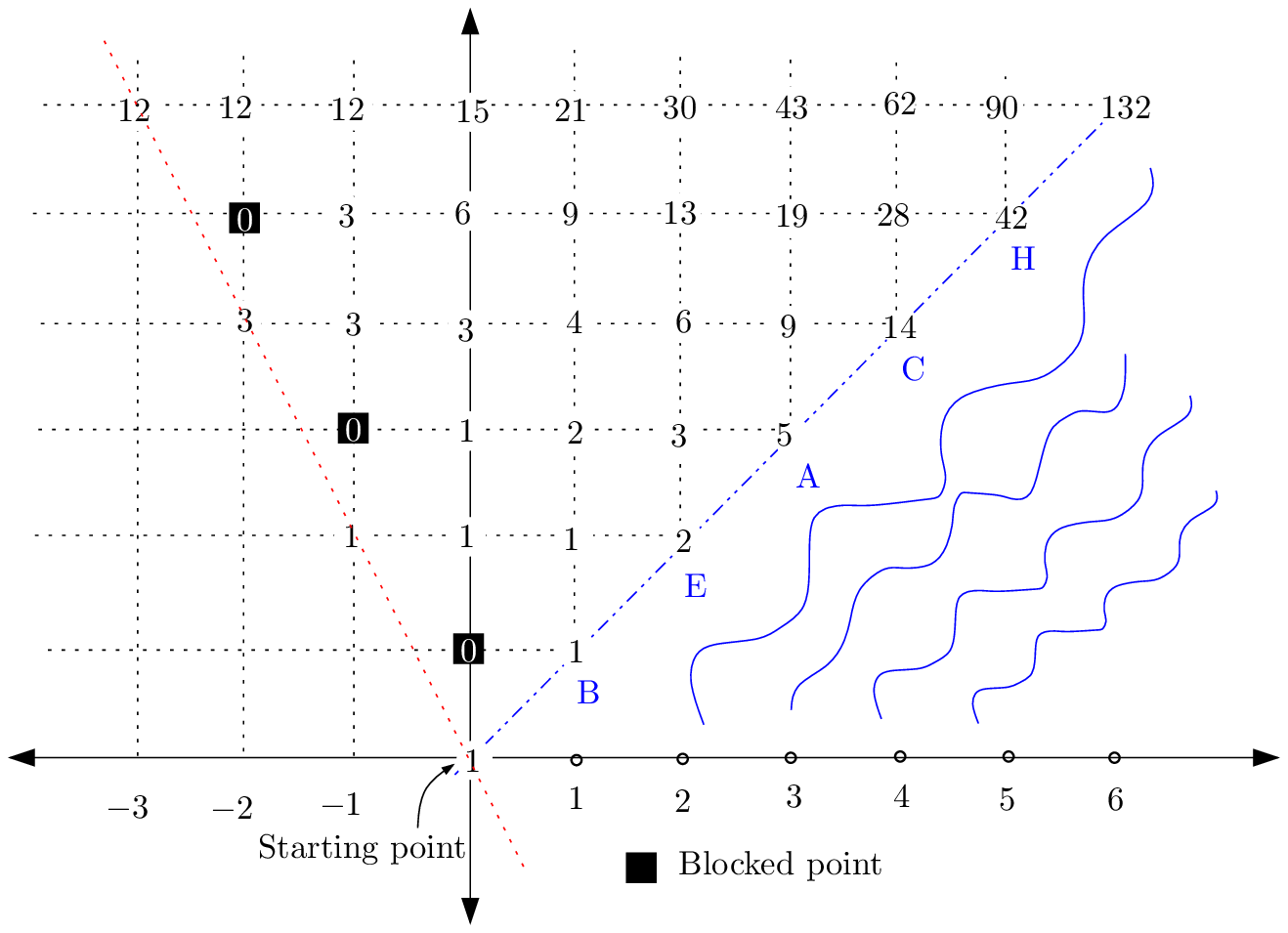}\\
 \medskip
 Fig.~2. Detours preserving Catalan traffic.
\end{center}

The main result of~\cite{Nied} is that the numbers along the ``beach"
are Catalan's number. This is proven in three different ways. One
of them relies on the following recursion:
\begin{eqnarray}
\left\{
\begin{array}{rll}
 \Upsilon(m,n)&=&\Upsilon(m+1,n)- \Upsilon(m,n-1)\\
 \Upsilon(n,n)&=&C_n \\
 \end{array}%
\right., \label{eq:rec}
\end{eqnarray} \noindent
holding in the domain $\{(m,n)\in\ZZ^2\,|\, -n\leq 2m\leq 2n\}$,
where  $\Upsilon(m,n)$  denotes the number of paths to get the
point $(m,n)$ starting from the origin. \claim{} The Catalan
number $C_n$ has also a beautiful geometric interpretation (see
e.g.~\cite{Mukai}): in fact, it is the Pl\"ucker degree
\[
\kappa_{2n,0}=\int_{G_1({\PP^{n+1}})}\sigma_1^{2n},
\]
of the  grassmannian of lines  $G_1(\PP^{n+1})$. The main result
of this paper is that {\em for each $(m,n)$, such that $n\geq
m\geq 0$}, $\Upsilon(m,n)=\kappa_{2m,n-m}$  where
\[
\kappa_{2m,n-m}=\int_{G_1({\PP^{n+1}})}\sigma_1^{2m}\sigma_2^{n-m},
\]
is the top intersection number computed in the integral cohomology
ring \linebreak
$H^*(G_1(\PP^{n+1}),\ZZ)$ of $G_1(\PP^{n+1})$, generated (as a
$\ZZ$-algebra) by the {\em special Schubert cycles} $\sigma_1$ and
$\sigma_2$. The proof consists in using the formalism introduced
in~\cite{Gat1} (see also~\cite{Gat2} and~\cite{Sant1}) to show
that the recursion~(\ref{eq:rec}) holds for $K(m,n):=\kappa_{2m,n-m}$.
\bigskip
\noindent {\bf Acknowledgment.} I want to thank my advisor
L.~Gatto for many helpful discussions and the anonymous referee
for ideas and comments.

\section{Preliminaries}\label{sec1}
\claim{}The grassmannian $G_1(\PP^{n+1})$ is a complete projective
variety of complex dimension $2n$. It is known (see
e.g.~\cite{Fu1}, \cite{GH}) that its {\em integral cohomology} (or
Chow intersection ring) $H^*(G_1(\PP^n),\ZZ)$ is generated by the
special Schubert cycles $\sigma_1$ and $\sigma_2$. The cycle
$\sigma_1\in H^*(G_1({\PP^{n+1}}),\ZZ)$ is the cohomology class
represented by the subvariety which is the closure of all the
lines incident a codimension $2$ linear subspace of $\PP^{n+1}$,
while $\sigma_2\in H^*(G_1({\PP^{n+1}}),\ZZ)$ is  represented by
the closure of all the lines which are incident a codimension $3$
linear subspace of $\PP^{n+1}$. A {\em top intersection number} in
$G_1(\PP^{n+1})$ is the {\em degree} of the product
$\sigma_1^a\sigma_2^b$ in $H^*(G_1({\PP^{n+1}}),\ZZ)$, with
$a+2b=2n$.

\claim{} The  cohomology (or the intersection) theory of complex
grassmannian varieties can be described by {\em Schubert
Calculus}. The latter can be phrased in purely algebraic terms
using the formalism introduced in~\cite{Gat1} (see \cite{Gat2} for
more details). For the grassmannian of lines of $\PP^{n+1}$ this
works, in short, as follows. Let $\bigwedge^2M$ be the $2^{nd}$
exterior power of a free module of rank $n+2$. If $M$ is spanned
by $(\ep^{n+1},\ep^n,\ldots,\ep^{1},\ep^0)$, then $\bigwedge^2M$
is freely generated by $\{\ep^i\wedge\ep^j\,|\, 0\leq i<j\leq
n+1\}$.

Let $D_1:M\lra M$ such that $D_1\ep^i=\ep^{i-1}$ if $i>1$ and
$D_{1}\varepsilon^{0}=0$. Extend $D_1$ to an endomorphism of
$\bigwedge^2M$ by setting:
\[
D_1(\ep^i\wedge\ep^j)=D_1\ep^i\wedge\ep^j+\ep^i\wedge D_1\ep^j,
\]
and let $D_2:\bigwedge^2M\lra \bigwedge^2M$ such that:
\[
D_2(\ep^i\wedge\ep^j)=D_1^2\ep^i\wedge\ep^j+D_1\ep^i\wedge D_1\ep^j+\ep^i\wedge D_1^2\ep^j.
\]
In other words $D_1$ behaves as a ``first derivative" and $D_2$ as
a ``second derivative". The main result in~\cite{Gat1} is that the
endomorphism $D_1$ and $D_2$ generate a commutative subalgebra
${\cal A}^*$ of the $\ZZ$-algebra $End_\ZZ(\bigwedge^2M)$ which is
isomorphic to $H^*(G_1(\PP^n),\ZZ)$. The isomorphism is explicitly
obtained from sending $ \sigma_1\mapsto D_1$ and $\sigma_2\mapsto
D_2$. From this point of view, it turns out that the degree of a
top intersection product $\sigma_1^a\sigma_2^b$ ($a+2b=2n$) in
$H^*(G_1(\PP^{n+1}),\ZZ)$ is nothing else than the coefficient
$\kappa_{a,b}$ in the equality:
\[
D_1^aD_2^b(\ep^{n+1}\wedge\ep^n)=\kappa_{a,b}\cdot \ep^1\w\ep^{0}.
\]
\section{The result}
In this section we will prove the main result of this paper: the
connection between the numbers in the Catalan Traffic and top
intersection numbers  in the integral cohomology ring of the
grassmannian $G_1(\PP^3)$ of lines in $\PP^3$. This connection is
a consequence of the following:

 \claim{\bf Theorem.}\label{teo} {\em Let
$K(m,n):=\kappa_{2m, n-m}$ be the coefficient of $\ep^1\w\ep^0$ in
the expansion of $D_{1}^{2m}D_{2}^{n-m}\left( \varepsilon
^{n+1}\wedge\varepsilon^{n}\right)$. Thus, for all $0\leq m\leq n$
the following recursion hold: \[K\left( m,n\right) =K\left(
m+1,n\right)  -K\left( m,n-1\right).\]
%
}

\proof Let $\Delta_{11}$ be the endomorphism of $M$ defined by:
\begin{equation} \Delta_{11}(D)(\ep^i\w\ep^j) = D_1\ep^i\w D_1\ep^j= (D_1^2
-D_2)(\ep^i\w\ep^j).\label{eq:ultt} \end{equation} Recalling that
\[
D_1^{2m}D_2^{n-m}(\ep^{n+1}\wedge\ep^{n})=K(m,n)\cdot\ep^{1}\wedge\ep^{0},
\]
by definition of $K(m,n)$, one has
\begin{eqnarray}
K(m,n)\cdot\ep^1\wedge\ep^{0}&=&D_1^{2m}D_2^{n-m}(\ep^{n+1}\wedge\ep^n)=\nonumber\\&=&D_1^{2m}D_2^{n-m-1}(D_1^2 -\Delta_{11})(\ep^{n+1}\wedge\ep^n)=\nonumber \\
                                    &=& D_1^{2m+2}D_2^{n-m-1}(\ep^{n+1}\wedge\ep^n)-D_1^{2m}D_2^{n-m-1}(\ep^{n}\wedge\ep^{n-1}).\label{eq:ult}
\end{eqnarray}
\noindent Now,  on the r.h.s of formula~(\ref{eq:ult}), the former
summand is precisely $K(m+1,n)\ep^1\wedge\ep^{0}$ while the
latter, using~(\ref{eq:ultt}),  is equal to
$D_1^{2m}D_2^{n-m-1}(\ep^n\wedge\ep^{n-1})$ which in turn equals
$\kappa_{2m,n-m-1}\ep^1\wedge\ep^{0}$. Hence, keeping in mind that
$\kappa_{2m,n-m-1}=K({m,n-1})$, one has, using~(\ref{eq:ult}):
\begin{eqnarray*}
K({m,n})\cdot\ep^{1}\wedge\ep^{0}
                                   = (K({m+1,n})-K({m,n-1}))\cdot
                                   \ep^{1}\wedge\ep^{0}.
\end{eqnarray*}
As a conclusion
$
K({m+1,n})=K({m,n})+K({m,n-1}).
$
\qed
\bigskip
\claim{\bf Proposition.}\label{prop}  {\em For all $n\geq0 $,
$K(n,n)=C_n$.}

\proof In fact one has:

\begin{eqnarray} D_1^{2n}(\ep^{n+1}\wedge\ep^{n})&=&\sum_{i=0}^{2n}{2n\choose
i}D_1^i\ep^{n+1}\wedge D_1^{2n-i}\ep^{n}\nonumber
\\&=&\sum_{i=0}^{2n}{2n\choose i}\ep^{n+1-i}\wedge
\ep^{i-n}\label{eq:exbin} \end{eqnarray}

Since $n\leq i \leq n+1$, $i=n$ or $i=n+1$. Hence
only the sum
\[
{2n\choose n} \ep^{1}\wedge \ep^0+{2n\choose n+1}\ep^{0}\wedge
\ep^{1}
\]
can survive in expression~(\ref{eq:exbin}). Therefore:
\begin{equation} D_1^{2n}\ep^{n+1}\wedge\ep^n=\left[{2n\choose n}-
{2n\choose n+1}\right]\ep^{1}\wedge \ep^{0} \end{equation} so that
\[
\hskip82ptK(n,n)={2n\choose n}- {2n\choose n+1}={(2n)!\over (n+1)!n!}=C_n.\hskip84pt\qed
\]

\bigskip

\claim{\bf Corollary.}\label{catath} {\em For all $0\leq m\leq n$,
the number $\Upsilon(m,n)$ (see equation (\ref{eq:rec})) coincides
with the number $\kappa_{2m, n-m}$ of lines in $\PP^{n+1}$
incident $2m$ linear subspaces  of codimension $2$ and $n-m$
subspaces of codimension $3$ in general position in $\PP^{n+1}$.}
\medskip

\proof This follows from Theorem \ref{teo},  Proposition
\ref{prop} and the remarks  in Section 2. \qed

\bigskip

\claim{}{\bf Remark.} A formula for the number $\kappa_{a,b}$,
with $a+2b=2n$ is computed in~\cite{Sant1}. For $a=2m$ and
$b=n-m$, it  coincides with the following combinatorial
expression:
\begin{eqnarray}
\kappa_{2m,n-m}&=& \sum_{i=0}^{n-m}\sum_{j=0}^{2m}{2m\choose
j}{(n-m)!(m+n-2j-3i+1) \over i!(n-j-2i+1)!( i+j-m)!},
\label{ftese}
\end{eqnarray}

\noindent holding for each $n\geq m\geq 0$.

Further manipulations involving other properties of the operators
$D_i$, studied in~\cite{Sant1} and not mentioned here, lead to the following simplified form:
\begin{eqnarray*}
\kappa_{2m,n-m}&=&\sum_{i=0}^{n-m}\left(-1\right) ^{i}{n-m\choose
i}C_{n-i},
\end{eqnarray*}

\noindent holding for each $n\geq m\geq 0$.

%
%

%

\bigskip

\begin{small}

\texttt{Dipartimento di Matematica, Politecnico di Torino\\
\indent  C.so Duca degli Abruzzi 24, 10129 - Torino-Italia \\
\indent E-mail address: \href{mailto: taise@calvino.polito.it}{
taise@calvino.polito.it}\\
\\
\indent Departamento de Ci\^{e}ncias Exatas e da Terra \\
\indent Universidade
Estadual de Feira de Santana \\
\indent Av. Universitária, s/n - Km 03 da BR 116 \\
\indent  44031-460,Feira de Santana - BA - Brasil \\ \indent
E-mail address: \href{mailto: taisesantiago@gmail.com}{
taisesantiago@gmail.com}}

\end{small}
\end{document}